\documentclass[12pt]{article}
\usepackage{amssymb,amsmath}
\usepackage{graphics}
\newtheorem{thm}{Theorem}
\newtheorem{cor}[thm]{Corollary}

\newenvironment{defin}{\medskip\noindent{\sc
Definition}.}{\goodbreak\medskip}
\newenvironment{nota}{\medskip\noindent{\sc
Notations}.}{\goodbreak\medskip}
\newenvironment{remk}{\noindent{\sc
Remark}.}{\goodbreak\vskip10pt}
\newtheorem{prop}[thm]{Proposition}

\newtheorem{ques}{Question}

\def\demo{\medskip\goodbreak\noindent
     \hbox{\sc Proof \kern .3em}\ignorespaces}%
  \def \qedbox{$\square$}%
  \def \qed{\hglue1mm\hfill{\ifmmode\qedbox
     \else\unskip\ \hglue0mm\hfill\qedbox\medskip
      \goodbreak\fi}}%
\def\enddemo{\qed\goodbreak\vskip10pt}%

\def\qed{\hglue1mm\hfill\raise -2pt\hbox{\vrule\vbox to 10pt{\hrule width
4pt
                  \vfill\hrule}\vrule}}

\newcommand{\T}{\mathbb {T}}

\newcommand{\esse}{\mathbb {S}}
\newcommand{\R}{\mathbb {R}}

\newcommand{\N}{\mathbb {N}}
\newcommand{\Nc}{\mathcal {N}}

\newcommand{\Cc}{\mathcal {C}}
\newcommand{\Ic}{\mathcal {I}}
\newcommand\calM{\mathfrak{M}}

\newcommand{\Mc}{\mathcal {M}}

\newcommand{\Gc}{\mathcal {G}}
\newcommand{\Lc}{\mathcal {L}}

\newcommand{\Ac}{\mathcal {A}}
\newcommand{\Tc}{\mathcal {T}}

\newcommand{\Sc}{\mathcal{S}}

\begin{document}
\title{{ A certain minimization property implies a certain integrability}}
\author{M.-C. ARNAUD
\thanks{ANR Project BLANC07-3\_187245, Hamilton-Jacobi and Weak KAM Theory}
\thanks{ANR DynNonHyp}
\thanks{Universit\'e d'Avignon et des Pays de Vaucluse, Laboratoire d'Analyse non lin\' eaire et G\' eom\' etrie (EA 2151),  F-84 018Avignon,
France. e-mail: Marie-Claude.Arnaud@univ-avignon.fr}
}

\maketitle
\abstract{The manifold $M$ being compact and connected and $H$ being a Tonelli Hamiltonian such that $T^*M$ is equal to the dual tiered Ma\~n\'e set, we prove that there is a partition of $T^*M$ into invariant $C^0$ Lagrangian graphs. Moreover, among these graphs, those that are $C^1$ cover a dense $G_\delta$ subset of $T^*M$. The dynamic restricted to each of these sets is non wandering.
}

\tableofcontents
\newpage
\section{Introduction}
In these article, we go on with our study of the so-called tiered Ma\~n\'e set. We began this study in \cite{Arna1}. Let us recall that the dual tiered Ma\~n\'e set $\Nc_*^T(H)$ of a Tonelli Hamiltonian\footnote{all these notions will be precisely defined in next section} is the union of all the dual Ma\~n\'e sets of $H$ associated to all the cohomology classes of $M$.\\
In \cite{Arna1}, we proved that for a generic Tonelli Hamiltonian, the tiered Ma\~n\'e set has no interior. \\
In our new article, we consider the following (non-generic) case~: we assume that ${\cal N}_*^T(H)=T^*M$.  In other words, we assume that every orbit of the Hamiltonian flow of $H$ is globally minimizing for $L-\lambda$, where $L$ is the Lagrangian associated to $H$ and $\lambda$ a closed 1-form (that depends on the considered orbit).\\
Such flows are part of a set of more general Tonelli Hamiltonian flows~: those that have no conjugate points. For example, it is proved in \cite{patpat1} that any Anosov Hamiltonian level of a Tonelli Hamiltonian has no conjugate points. The same result for geodesic flows was proved in the 70's  by W.~Klingenberg in \cite{Klin1}. But the tiered Ma\~n\'e set of an Anosov geodesic flow has no interior (see \cite{Arna1}) hence in this case, the dual tiered Ma\~n\'e set is not equal to $T^*M$. In fact, we  prove~:

\begin{thm}\label{T1}
Let $M$ be a compact and connected manifold and let $H~: T^*M\rightarrow \R$ be a Tonelli Hamiltonian. Then the two following assertions are equivalent~:
\begin{enumerate}
\item there exists a partition of $T^*M$ into  invariant Lipschitz Lagrangian graphs;
\item the dual tiered Ma\~n\'e set of $H$ is the whole cotangent bundle $T^*M$.
\end{enumerate}
Moreover, in this case~:
\begin{enumerate}
\item[$\bullet$] there exists an invariant dense $G_\delta$ subset $\Gc$ of $T^*M$ such that all the graphs of the partition that meets $\Gc$ are in fact $C^1$.
\item[$\bullet$]  Mather's $\beta$ function is everywhere differentiable.
\end{enumerate}
\end{thm}
Let us emphasize why this result is surprizing~: we just ask that all the orbits are, in a certain way, minimizing, and we prove that they are well-distributed on invariant Lipschitz Lagrangian graphs.

An easy corollary is the following~:
\begin{cor}\label{C2}
Let $M$ be a compact and connected manifold and let $H~: T^*M\rightarrow \R$ be a Tonelli Hamiltonian.  Then the two following assertions are equivalent~:
\begin{enumerate}
\item there exists a partition of $T^*M$ into invariant Lipschitz  Lagrangian graphs;
\item  $T^*M$ is covered by the union of its  invariant Lipschitz Lagrangian graphs.
\end{enumerate}
The same statement is true if we replace everywhere ``Lipschitz'' by ``smooth''.
\end{cor}
In \cite{Arna2}, we proved a Birkhoff multimensional theorem for Tonelli Hamiltonians. We deduce~:

\begin{cor}\label{C3}
Let $M$ be a closed and connected manifold and let $H~: T^*M\rightarrow \R$ be a Tonelli Hamiltonian.  Then the two following assertions are equivalent~:
\begin{enumerate}
\item  there exists a partition of $T^*M$ into Lagrangian invariant smooth graphs;
\item $T^*M$ is covered by the union of its Lagrangian invariant smooth submanifolds that are Hamiltonianly isotopic to some Lagrangian smooth graph.
\end{enumerate}
\end{cor}
These results give us a characterization of a weak form of integrability; following \cite{Arna3}, we say that a Tonelli Hamiltonian is {\em $C^0$-integrable} if there is a partition of $T^*M$ into invariant $C^0$-Lagrangian graphs, one for each cohomology 
class in $H^1(M, \R)$. We then prove that if all the orbits are in some Ma\~n\'e set, then the Hamiltonian is $C^0$-integrable. A natural question is then~:
\begin{ques}~: does there exist any Tonelli Hamiltonian that is $C^0$-integrable but not $C^1$-integrable (i.e. for which the invariant graphs are not all $C^1$)?
\end{ques}

Let us notice that we finally prove that our hypotheses implies that the function $\beta$ is everywhere differentiable. An interesting question, well-known from   specialists, is~: when the function $\beta$ is everywhere differentiable, is the Hamiltonian $C^0$-integrable? In the case of closed surfaces, a positive answer to this question is given in \cite{masor}.

\noindent{\em Part of this work was done at the University of Maryland in April 2010}
\medskip 

\section{An overview of Mather-Ma\~n\'e-Fathi theory of minimizing orbits}
 \subsection{Tonelli Lagrangian and Hamiltonian functions}
 Let $M$ be a compact and  connected  manifold endowed with a
Riemannian metric . We denote a point of the tangent bundle $TM$ by
$(q, v)$  with $q\in M$ and $v$ a vector tangent to $M$ 
at
$q$. The projection 
$\pi: TM\rightarrow M$ is then $(q, v)\rightarrow q$. The notation $(q, p)$  
designates a point of the cotangent bundle $T^*M$ with $p\in T^*_qM$  and $\pi^*: T^*M\rightarrow M$ is the canonical projection $(q , p)\rightarrow q$.

We consider a Lagrangian  function $L: TM\rightarrow \R$ which is $C^2$
and: \begin{enumerate}
\item[$\bullet$] uniformly superlinear: uniformly on $q\in M$, we have: $\displaystyle{\lim_{\| v\|\rightarrow +\infty}\frac{L(q, v)}{\| v\|} =+\infty}$;
\item[$\bullet$] strictly convex: for all $(q, v)\in TM$, $\frac{\partial^2
L}{\partial v^2}(q, v)$ is positive definite.
\end{enumerate} 
Such a Lagrangian function will be called a {\em Tonelli Lagrangian
function}.\\

We can associate to such a Lagrangian  function  the Legendre map
$\Lc=\Lc_L: TM\rightarrow T^*M$ defined by: $\Lc(q, v)=\frac{\partial
L}{\partial v}(q, v)$ which is a fibered
$C^2$ diffeomorphism and the    Hamiltonian function $H: T^*M\rightarrow
\R$ defined by: $H(q,p)=p\left( \Lc^{-1}(q,p)\right) -L(\Lc^{-1}(q,p))$
(such a Hamiltonian function will be called a {\em Tonelli Hamiltonian
function}). The Hamiltonian function 
$H$ is then superlinear, strictly convex in the fiber and $C^2$. We
denote by
$(f^L_t)$ or $(f_t)$ the Euler-Lagrange flow associated to $L$ and $(\varphi_t^H)$ or $(\varphi_t)$  the
Hamiltonian flow associated to $H$; then we have~: $\varphi_t^H=\Lc\circ
f_t^L\circ\Lc^{-1}$.\\

If $\lambda$ is a ($C^\infty$) closed 1-form of $M$, then the map
$T_\lambda~: T^*M\rightarrow T^*M$ defined by~: $T_\lambda (q, p)=(q,
p+\lambda (q))$ is a symplectic ($C^\infty$) diffeomorphism; therefore,
we have~: $(\varphi^{H\circ T_\lambda}_t)=(T_\lambda^{-1}\circ \varphi_t\circ
T_\lambda )$, i.e.  the Hamiltonian flow of $H$ and $H\circ T_\lambda$ are
conjugated. Moreover, the Tonelli Hamiltonian function $H\circ T_\lambda$ is
associated to the Tonelli Lagrangian function $L-\lambda$, and it is
well-known that~: $(f_t^L)=(f_t^{L-\lambda})$; the two Euler-Lagrange
flows are equal. Let us emphasize that these flows are equal, but the
Lagrangian functions, and then the Lagrangian actions differ and so  the
minimizing ``objects'' may be different.
\subsection{Tiered sets~: Mather, Aubry and Ma\~n\'e}\label{sect22}
For   a Tonelli Lagrangian function ($L$ or $L-\lambda$), J.~Mather introduced in \cite{Mat3}
(see  \cite{M2} too) a particular subset $\Ac(L-\lambda)$ of $TM$
which he called the ``static set'' and which is now usually called the
``{\em Aubry set}'' (this name is due to A.~Fathi)\footnote{ These sets extend the notion of
``Aubry-Mather'' sets for the twist maps.}. There exist  different but equivalent definitions of
this set (see
\cite{CIPP} ,  
\cite{Fa1},
\cite{M2} and subsection \ref{sect23}) and it is known that two closed 1-forms
which are in  the same cohomological class define the same Aubry set~:
$$[\lambda_1]=[\lambda_2]\in H^1(M)\Rightarrow \Ac (L-\lambda_1)=\Ac
(L-\lambda_2).$$
We can then   introduce the following notation~: if $c\in H^1(M)$ is
a cohomological class, 
$\Ac_c(L)=\Ac (L-\lambda)$ where $\lambda$ is any closed 1-form belonging
to $c$. $A_c(L)$ is compact, non empty and invariant under $(f_t^L)$.
Moreover, J.~Mather proved in \cite{Mat3} that it is  a Lipschitz graph  above a part of the
zero-section (see
\cite{Fa1} or subsection \ref{sect23} too).

As we are as interested in the Hamiltonian dynamics as well as in the
Lagrangian ones, let us  define the dual Aubry set~:

\begin{enumerate}
\item[--] if $H$ is the Hamiltonian function associated to the Tonelli
Lagrangian function $L$, its {\em dual Aubry set} is $\Ac^*(H)=\Lc_L
(\Ac (L))$; 
\item[--] if $c\in H^1(M)$ is a cohomological class, then
$\Ac^*_c(H)=\Lc_L(\Ac_c(L))$ is the {\em $c$-dual Aubry set}; let us
notice that for any closed 1-form
$\lambda$ belonging to $c$, we have~: $T_\lambda(
\Ac^*(H\circ T_\lambda))=\Ac_c^*(H)$.
\end{enumerate}
These sets are invariant under the Hamiltonian flow $(\varphi_t^H)$. 

Another important invariant subset in the theory of Tonelli Lagrangian
functions is the so-called Mather set.  For it, there exists one
definition (which is in     \cite{Fa1},
\cite{M2}, \cite{mather1} and subsection \ref{sect24})~: it is the closure of the union of the supports
of the minimizing measures for $L$; it is denoted by $\Mc(L)$ and the {\em
dual Mather set} is $\Mc^*(H)=\Lc_L(\Mc (L))$ which is compact,
non empty and invariant under the flow $(\varphi_t^H)$. As for the Aubry set,
if
$c\in H^1(M)$ is a cohomological class, we define~: $\Mc_c(L)=\Mc
(L-\lambda)$ which is independent of the choice of the closed 1-form
$\lambda$ belonging to $c$. Then
$\Mc^*_c(H)=\Lc_L(\Mc_c(L))=T_\lambda(\Mc^*(H\circ T_\lambda))$ is
invariant under $(\varphi^H_t)$; we name it the {\em $c$-dual Mather set}. 

In a similar way, if  $\Nc (L)$ is  the Ma\~n\'e set, the {\em dual
Ma\~n\'e set} is 
$\Nc^*(H)=\Lc_L(\Nc(L))$; we note that if $c\in H^1(M)$ and $\lambda\in
c$, then $\Nc_c(L)=\Nc(L-\lambda)$ is independent  of the choice of
$\lambda\in c$ and then the {\em $c$-dual Ma\~n\'e set} is
$\Nc_c^*(H)=\Lc_L(\Nc_c(L))=T_\lambda(\Nc^*(H\circ T_\lambda))$; it  is
invariant under $(\varphi_t^H)$, compact and non empty but is not
necessarily a graph.

For every cohomological class $c\in H^1(M)$, we have the inclusion~:
$\Mc^*_c(H)\subset \Ac^*_c(H)\subset  \Nc^*_c(H) $. Moreover, there
exists a real number denoted by
$\alpha_H (c)$ such that~: $\Nc^*_c(H)\subset H^{-1} (\alpha_H (c))$ (see
\cite{Car} and \cite{mather1}), i.e. each dual Ma\~n\'e set is contained in an energy level.
For $c=0$, the value $\alpha_H(0)$ is   named
the ``critical value'' of $L$.\\

\begin{defin} If $H~: T^*M\rightarrow \R$ is a Tonelli Hamiltonian
function, the {\em tiered Aubry set}, the {\em tiered Mather set}  and
the {\em tiered Ma\~n\'e set} are~:
$$\Ac^T(L)=\bigcup_{c\in H^1(M)}\Ac_c(L);\quad \Mc^T(L)=\bigcup_{c\in
H^1(M)}\Mc_c(L);\quad \Nc^T(L)=\bigcup_{c\in
H^1(M)}\Nc_c(L).$$
Their dual sets are~:
$$\Ac^T_*(H)=\bigcup_{c\in H^1(M)}\Ac^*_c(H);\quad
\Mc^T_*(H)=\bigcup_{c\in H^1(M)}\Mc^*_c(H);\quad \Nc^T_*(H)=\bigcup_{c\in
H^1(M)}\Nc^*_c(H).$$
\end{defin}

\subsection{Ma\~n\'e potential, Peierls barrier, Aubry and Ma\~n\'e sets}\label{sect23}
We gather in this sections some well-known results; the ones concerning
the Peierls barrier are essentially due to A.~Fathi (see
\cite{Fa1}), the others concerning Ma\~n\'e potential are given in
\cite{M1}, \cite{CDI} and \cite{CI}.\\
  In the whole section,
$L$ is a Tonelli Lagrangian function.

\begin{nota}\begin{enumerate}
\item[$\bullet$] given two points $x$ and $y$ in $M$ and $T>0$, we
denote by $\Cc_T(x,y)$ the set of absolutely continuous curves $\gamma~:
[0, T]\rightarrow M$ with $\gamma (0)=x$ and $\gamma (T)=y$;
\item[$\bullet$]  the Lagrangian action along an absolutely continuous
curve
$\gamma~: [a, b]\rightarrow M$ is defined by~:
$$A_L(\gamma)=\int_a^bL(\gamma(t), \dot\gamma (t))dt;$$
\item[$\bullet$] for each $t>0$, we define the function $h_t~:
M\times M\rightarrow \R$ by~: 
$h_t(x,y)=\inf\{ A_{L+\alpha_H(0)}(\gamma);
\gamma\in
 \Cc_t(x,y)\}$;
\item[$\bullet$] the Peierls barrier is then the function $h~: M\times
M\rightarrow \R$ defined by~:
$$h(x,y)=\liminf_{t\rightarrow + \infty} h_t(x,y);$$ 
\item[$\bullet$]   we define the {\em (Ma\~n\'e)  
potential} $m~: M\times M\rightarrow \R$ by~: $m 
(x,y)=\inf\{ A_{L+\alpha_H(0)}(\gamma);
\gamma\in
\bigcup_{T>0}\Cc_T(x,y)\}=\inf\{ h_t(x,y); t>0\}$.
\end{enumerate}
\end{nota}
Then, the Ma\~n\'e potential verifies~: 
\begin{prop} We have~:
\begin{enumerate}
\item $m$ is finite and $m\leq h$;
\item $\forall x, y, z\in M, m(x,z)\leq m(x,y)+m(y,z)$;
\item $\forall x\in M, m(x,x)=0$;
\item if $x,y\in M$, then $m(x,y)+m(y,x)\geq 0$;
\item if $M_1=\sup\{ L(x,v); \| v\| \leq 1\}$, then~: $\forall x, y\in
M, |m(x,y)|\leq (M_1+\alpha_H(0))d(x,y)$;
\item $m~: M\times M\rightarrow \R$ is $(M_1+\alpha_H(0))$-Lipschitz. 
\end{enumerate}
\end{prop} 
Now we can define~: 

\begin{defin}\begin{enumerate}
\item[$\bullet$] a absolutely continuous curve $\gamma~:
I\rightarrow M$ defined on an interval $I$ is a {\em ray} if~:
$$\forall [a, b]\subset I,
A_{L+\alpha_H(0)}(\gamma_{|[a,b]})=h_{(b-a)}(\gamma(a), \gamma(b));$$
a ray is always a solution
of the Euler-Lagrange equations;
\item[$\bullet$] a absolutely continuous curve $\gamma~:
I\rightarrow M$ defined on an interval $I$ is {\em semistatic} if~:
$$\forall [a, b]\subset I, m(\gamma (a), \gamma
(b))=A_{L+\alpha_H(0)}(\gamma_{|[a,b]});$$ a semistatic curve is always a ray;
\item[$\bullet$] the {\em Ma\~n\'e set} is then~: $\Nc(L)=\{ v\in TM;
\gamma_v\quad is \quad semistatic\}$ where $\gamma_v$ designates the
solution $\gamma_v~: \R\rightarrow M$  of the Euler-Lagrange equations with initial condition  $v$
for
$t=0$; $\Nc (L)$  is  contained in the critical energy level;
\item[$\bullet$] a absolutely continuous curve $\gamma~:
I\rightarrow M$ defined on an interval $I$ is {\em  static} if~:
$$\forall [a, b]\subset I, -m(\gamma (b), \gamma
(a))=A_{L+\alpha_H(0)}(\gamma_{|[a,b]});$$ a  static curve is always a semistatic
curve;
\item[$\bullet$] the {\em Aubry set} is then~: $\Ac(L)=\{ v\in TM;
\gamma_v\quad is \quad  static\}$.

\end{enumerate}
\end{defin}
The following result is proved in \cite{CI}~:
\begin{prop}\label{Pr4}
If $v\in TM$ is such that $\gamma_{v|[a, b]}$ is static for some $a<b$,
then $\gamma_v~: \R\rightarrow M$ is static, i.e. $v\in \Ac(L)$.
\end{prop}\label{Peierls}
The Peierls barrier verifies (this proposition contains some
results of \cite{CP}, \cite{Fa1} and \cite{Co1})~:
\begin{prop} (properties of the Peierls barrier $h$) \label{PPeierls}
\begin{enumerate}
\item the values of the map $h$ are finite  and $m\leq h$;
\item if $M_1=\sup\{ L(x,v); \| v\| \leq 1\}$, then~:
$$\forall x,y,x',y'\in M, |h(x,y)-h(x',y')|\leq
(M_1+\alpha_H(0))(d(x,x')+d(y,y'));$$   
therefore $h$ is Lipschitz;
\item  if $x,y\in M$, then $h(x,y)+h(y,x)\geq 0$; we deduce~: $\forall
x\in M, h(x,x)\geq 0$;
\item $\forall x, y, z\in M, h(x,z)\leq h(x,y)+h(y,z)$;
\item \label{I5} $\forall x\in M, \forall y\in \pi(\Ac (L)),
m(x,y)=h(x,y)\quad{\rm and}\quad m(y,x)=h(y,x)$;
\item $\forall x\in M,h(x,x)=0\Longleftrightarrow x\in \pi(\Ac(L))$.
\end{enumerate}
\end{prop}  
 The last item of this
proposition gives us a characterization of the projected Aubry set
$\pi(\Ac(L))$. Moreover,  we have~:
\begin{prop}\label{Pr7} (A.~Fathi, \cite{Fa1}, 6.3.3) When $t$ tends to
$+\infty$, uniformly on $M\times M$, the function $h_t$ tends to the
Peierls barrier
$h$.
\end{prop}
A corollary of this result is given in \cite{CI}~:
\begin{cor} (\cite{CI}, 4-10.9) All the rays defined on $\R$ are semistatic.

\end{cor}

Let us give some properties of the Aubry and Ma\~n\'e sets (see \cite
{M2} and \cite{CDI})~:

\begin{prop}\label{Pr8} Let $L~: TM\rightarrow \R$ be a Tonelli Lagrangian
function. Then~:
\begin{enumerate}
\item[$\bullet$] the Aubry and Ma\~n\'e set are compact, non empty  and
$\Ac(L)\subset
\Nc (L)$;
\item[$\bullet$] the Aubry set is a Lipschitz graph above a part of the
zero section;
\item[$\bullet$] if $\gamma~: \R\rightarrow M$ is semistatic, then
$(\gamma, \dot\gamma)$ is a Lipschitz graph above a part of the
zero section;
\item[$\bullet$] the $\omega$ and $\alpha$-limit sets of every point of the Ma\~n\'e set are
contained in the Aubry set.
\end{enumerate}

\end{prop}

Last item in proposition \ref{PPeierls} gives us a criterion to some $q\in M$ belong to some projected Aubry set. We will need a little more than this~: we will need to know what happens for its lift, the   Aubry set.

\begin{prop}\label{Pradial}  Let $c\in H^1(M)$ and $\lambda\in c$, $\varepsilon >0$  and let
$L~: TM\rightarrow
\R
$ be a Tonelli Lagrangian function. Then there exists $T_0>0$ such
that~: \\
$\forall T\geq T_0, \forall (q_0,v_0)\in \Ac_c(L,), \forall \gamma~: [0, T]\rightarrow M$ minimizing for $L-\lambda$ between $q_0$ and $q_0$, i.e.~:\\
 $\forall \eta~: [0, T]\rightarrow M, \eta (0)=\eta (T)=q_0\Rightarrow  \int_0^T(L(\gamma, \dot\gamma)-\lambda (\dot\gamma)+\alpha_H(c))\leq  \int_0^T(L(\eta, \dot\eta)-\lambda (\dot\eta)+\alpha_H(c))$\\
 then we have~: $d((q_0, v_0), (q_0, \gamma'(0)))\leq \varepsilon$

\end{prop}
\demo
Let us assume that the result is not true; then we may find a sequence
$(T_n)_{n\in \N}$ in $\R^*_+$ tending to $+\infty$, a  sequence $\gamma_n~: [0, T_n]\rightarrow M$ of absolutely continuous loops, all of whose minimizing   for $L-\lambda$ from $q_n$ to $q_n$ where $(q_n, w_n)\in \Ac_c(L)$ such that the sequence $(q_n, v_n)=(\gamma_n(0), \dot\gamma_n(0))$ satifies~: $\forall n\in\N, d((q_n, v_n), (q_n, w_n))\geq \varepsilon$.\\
The sequence $(q_n, v_n)$ is bounded (it is a consequence of the
so-called ``a priori compactness lemma'' (see \cite{Fa1}, corollary
4.3.2)); therefore we may extract a converging subsequence; we call it
$(q_n, v_n)$ again and $(q_\infty, v_\infty)$ is its limit.  Then $q_\infty\in \pi(\Ac_c(L))$ because the Aubry set is closed. We denote by $(q_\infty, w_\infty)\in \Ac_c(L)$ its lift. Then $\displaystyle{w_\infty=\lim_{n\rightarrow \infty}w_n}$ because $\Ac_c(L)$ is closed. Then~: $d((q_\infty, v_\infty),
(q_\infty, w_\infty))\geq \varepsilon$. \\
Now we use   proposition \ref{Pr7}~: we know that if we define $h_t^\lambda~:
M\times M\rightarrow \R$ by $\displaystyle{h^\lambda_t(x,y) =\inf\{
A_{L-\lambda+\alpha_H(c)}(\gamma);
\gamma\in
 \Cc_t(x,y)\}}$ and $\displaystyle{h^\lambda (x,y)=\liminf_{t\rightarrow +\infty}
h^\lambda _t(x,y)}$, the functions $h^\lambda_t$ tend    uniformly  to
$h^\lambda$ when $t$ tends to $+\infty$;  we have then~: \\
$ h_{T_n}^\lambda
(q_n,q_n)=A_{L-\lambda +\alpha_H(c)}(\gamma_n)$ tends to $h^{\lambda}(q_\infty, q_\infty)=0$ when $n$
tends to the infinite.\\
Let $\gamma_\infty$ be the solution of the
Euler-lagrange equations such that $(\gamma_\infty(0),
\dot\gamma_\infty(0))=(q_\infty, v_\infty)$. We want to prove that
$\gamma_\infty$ is static~: we shall obtain a contradiction. When $n$ is
big enough,
$\gamma_n(T_n)=\gamma_n(0)$  is close to $q_\infty$ and $\gamma_n(1)$ is close to $\gamma_\infty (1)$.
Let us fix $\eta>0$; then we define $\Gamma_{n}^{\eta}~: [0,
T_n+2\eta ]\rightarrow M$ by~: 
\begin{enumerate}
\item [$\bullet$] $\Gamma^ \eta_{n |[0, 1]}=\gamma_{\infty|[0,
1]}$;
\item [$\bullet$]  $\Gamma^ \eta_{n |[1, 1+\eta]}$ is a short
geodesic joining $\gamma_\infty (1)$ to $\gamma_n(1)$;
\item[$\bullet$] $\forall t\in [1+\eta, T_n+\eta], \Gamma^\eta_n
(t)=\gamma_n(t-\eta)$;
\item[$\bullet$] $\Gamma^ \eta_{n |[T_n+\eta, T_n+2\eta]}$ is a short
geodesic joining $\gamma_n(T_n)$ to $\gamma_\infty (0)$.
\end{enumerate}
If we choose carefully a sequence $(\eta_n)$ tending to $0$, we have~: 
$$\lim_{n\rightarrow \infty} A_{L-\lambda +
\alpha_H(c)}(\Gamma_n^{\eta_n})=\lim_{n\rightarrow \infty} A_{L-\lambda +
\alpha_H(c)}(\gamma_n)=0.$$
Because the contribution to the action of the two small geodesic arcs
tends to zero (if the $\eta_n$ are well chosen), this implies~:
$$A_{L-\lambda+\alpha_H(c)}(\gamma_{\infty |[0,1]})+m^\lambda(\gamma_\infty (1),
\gamma_\infty (0))\leq 0,$$ where $m^\lambda$ designates Ma\~n\'e potential
for the Lagrangian function $L-\lambda$. We deduce then from the
definition of Ma\~n\'e potential that $m^\lambda (\gamma_\infty (0),
\gamma_\infty (1))+ m^\lambda(\gamma_\infty (1),
\gamma_\infty (0))= 0$ and that~:
$A_{L-\lambda+\alpha_H(c)}(\gamma_{\infty|[0,1]})=m^\lambda (\gamma_\infty (0),
\gamma_\infty (1))$. It implies then that
$A_{L-\lambda+\alpha_H(c)}(\gamma_{\infty |[0,1]})=- m^\lambda(\gamma_\infty (1),
\gamma_\infty (0))$. Let us notice that, changing slightly
$\Gamma_n^\eta$, we obtain too~: 
$$\forall [a, b]\subset [0, +\infty[,
A_{L-\lambda+\alpha_H(c)}(\gamma_{\infty|[a,b]})=- m^\lambda(\gamma_\infty
(b),
\gamma_\infty (a));$$
therefore $\gamma_{\infty|[0, +\infty[}$ is static. To conclude, we use 
proposition \ref{Pr4}. \enddemo

\subsection{Minimizing measures, Mather $\alpha$ and $\beta$ functions}\label{sect24} The general references for this section are \cite{mather1} and \cite{masor}.
Let $\calM (L)$ be the space of compactly supported Borel probability measures invariant under the Euler-Lagrange flow $(f_t^L)$.  To every $\mu\in \calM (L)$ we may associate its average action $A_L(\mu)=\int_{TM}Ld\mu$. It is proved in \cite{mather1} that for every $f\in C^1(M, \R)$, we have~: $\int df(q).v d\mu(q,v)=0$. Therefore we can define on $H^1(M, \R)$ a linear functional $\ell(\mu)$ by~: $\ell(\mu)([\lambda])=\int\lambda (q).vd\mu (q,v)$ (here $\lambda$ designates any closed 1-form). Then  there exists a unique element $\rho(\mu)\in H_1(M, \R)$ such that~:
$$\forall \lambda, \int_{TM}\lambda (q).vd\mu(q,v)=[\lambda].\rho(\mu).$$
The homology class  $\rho(\mu)$ is called the {\em rotation vector} of $\mu$. Then the map $\mu\in\calM (L)\rightarrow \rho(\mu)\in H^1(M, \R)$ is   onto.  We can then define Mather $\beta$-function $\beta~: H_1(M, \R)\rightarrow \R$ that associates  the minimal value of the average action $A_L$ over the set of measures of $\calM(L)$ with rotation vector $h$ to each homology class $h\in H_1(M, \R)$. We have~:
$$\beta (h)=\min_{\mu\in\calM (L); \rho(\mu)=h}A_L(\mu).$$
A measure $\mu\in \calM (L)$ realizing such a minimum, i.e. such that $A_L(\mu)=\beta (\rho(\mu))$ is called a {\sl minimizing measure with rotation vector} $\rho (\mu)$. The $\beta$ function is convex and superlinear,  and we can define its conjugate function (given by Fenchel duality) $\alpha~: H^1(M, \R)\rightarrow \R$  by~:
$$\alpha ([\lambda]) =\max_{h\in H_1(M, \R)}([\lambda].h-\beta (h))=-\min_{\mu\in\calM(L)}A_{L-\lambda}(\mu).$$
A measure $\mu\in\calM (L)$ realizing the minimum of $A_{L-\lambda }$ is called a {\sl $[\lambda]$-minimizing measure}.\\
Being convex, Mather's $\beta$ function has a subderivative at any point $h\in H_1(M, \R)$; i.e. there exists $c\in H^1(M, \R)$ such that~: $\forall k\in H_1(M, \R), \beta (h)+c.(k-h)\leq \beta (k)$. We denote by $\partial \beta (h)$ the set of all the subderivatives of $\beta$ at $h$. By Fenchel duality, we have~:
$c\in\partial\beta (h)\Leftrightarrow c.h=\alpha (c)+\beta (h)$.\\
Then we introduce the following notations~:
\begin{enumerate}
\item[$\bullet$] if $h\in H_1(M, \R)$, the Mather set for the rotation vector $h$ is~: $$\displaystyle{\Mc^h(L)=\bigcup\{ {\rm supp}\mu; \quad\mu \quad{\rm is}\quad{\rm minimizing}\quad{\rm with}\quad{\rm rotation}\quad{\rm vector}\quad h\}};$$
\item[$\bullet$] if $c\in H^1(M, \R)$, the Mather set for the cohomology class $c$ is~: 
$$\displaystyle{\Mc_c(L)=\bigcup\{ {\rm supp}\mu; \quad\mu \quad{\rm is}\quad c-{\rm minimizing}\}}.$$
\end{enumerate}
The following equivalences are proved in \cite{masor} for any pair $(h, c)\in H_1(M, \R)\times H^1(M, \R)$~: 
$$\Mc^h(L)\cap \Mc_c(L)\not=\emptyset\Leftrightarrow \Mc^h(L)\subset \Mc_c(L)\Leftrightarrow c\in\partial\beta (h).$$
As explained in subsection \ref{sect22}, the dual Mather set for the cohomology class $c$ is defined by~: 
$\Mc_c^*(H)=\Lc_L(\Mc_c(L))$. If $\cal M^*(H)$ designates the set of compactly supported Borel probability measures of $T^*M$ that are invariant by the Hamiltonian flow $(\varphi_t)$, then the map $\Lc_*~: \calM(L)\rightarrow \calM^*(H)$ that push forward the measures by $\Lc$ is a bijection. We denote $\Lc_*(\mu)$ by  $\mu^*$ and say that the measures are dual. We say too that $\mu^*$ is minimizing if $\mu$ is minimizing in the previous sense. \\
Moreover, the Mather set $\Mc_c^*(H)$ is a subset of the Ma\~n\'e set $\Nc_c^*(H)$ and every invariant Borel probability measure the support of whose is in $\Nc_c^*(H)$ is $c$-minimizing.
\subsection{The link with the weak KAM theory}\label{secweak}
If $\lambda$ is a closed 1-form on $M$, we can consider the Lax-Oleinik semi-groups of $L-\lambda$, defined on $C^0(M, \R)$ by~:
\begin{enumerate}
\item[$\bullet$] the negative one~: $T_t^{\lambda, -}u=\min \left(u(\gamma(0))+\int_0^t(L(\gamma(s), \dot\gamma(s))-\lambda(\gamma(s))\dot\gamma (s))ds\right);$
where the infimum is taken on the set of $C^1$ curves  $\gamma~: [0, t]\rightarrow M$ such that   $\gamma (t)=q $;
\item[$\bullet$] the positive one~:
$  T^{\lambda,+}_tu(q)=\max  \left(u(\gamma(t))-\int_0^t(L(\gamma (s), \dot\gamma(s)-\lambda (\gamma (s)).\dot\gamma (s)))ds\right);$
where the infimum is taken on the set of $C^1$ curves  $\gamma~: [0, t]\rightarrow M$ such that $\gamma (0)=q $. 
\end{enumerate}
A.~Fathi proved in \cite{Fa1} that for each closed 1-form $\lambda$, there exists $k\in\R$ and $u\in C^0(M, \R)$ such that~: $\forall t>0, T_t^{\lambda, -}u=u-kt$ (resp. $\forall t>0, T_t^{\lambda, +}u=u+kt$). In this case, we have~: $k=\alpha ([\lambda])$. The function $u$ is called a negative (resp. positive) weak KAM solution for $L-\lambda$. We denote the set of negative (resp. positive) weak KAM solutions for $L-\lambda$ by $\Sc^-_\lambda$ (resp. $\Sc^+_\lambda$). \\
 Moreover, it is proved too that a function $u~: M\rightarrow \R$ that is $C^{1}$ is a positive weak KAM solution if and only if it is a negative weak KAM solution if and only if it is a   solution of the Hamilton-Jacobi equation~: $H(q, \lambda (q)+du(q))=\alpha ([\lambda])$. It is equivalent too to the fact that the graph of $\lambda+du$ is invariant by the Hamiltonian flow $(\varphi_t^H)$.\\
 But in  general, the weak KAM solutions are not $C^1$ and the graph of $\lambda+du$ is not invariant by the Hamiltonian flow. There is an invariant subset contained in all these graphs~: the dual Aubry set. Let us now recall which characterization of this set is given by A.~Fathi in \cite{Fa1}.\\
 A pair $(u_-, u_+)$ of negative-positive weak KAM solution is called a pair of conjugate weak KAM solutions if $u_{-|\pi(\Mc(L-\lambda))}=u_{+|\pi(\Mc(L-\lambda))}$. Each negative weak KAM solution has an unique conjugate positive weak KAM solution, and we define for any pair $(u_-, u_+)\in \Sc^-_\lambda\times\Sc^+_\lambda$ of conjugate weak KAM solutions for $L-\lambda$~:
 \begin{enumerate}
 \item[$\bullet$] $\Ic(u_-, u_+)=\{ q\in M, u_-(q)-=u_+(q)\}$;
 \item[$\bullet$] $\tilde\Ic(u_-, u_+)=\{ (q, du_-(q)); q\in \Ic(u_-, u_+)\}= \{ (q, du_+(q)); q\in \Ic(u_-, u_+)\}$.
 \end{enumerate}
 Then~: $\Ac^*_{[\lambda]}(H)=T_\lambda( \bigcap \tilde\Ic(u_-, u_+))$ where the intersection is taken on all the pairs of conjugate weak KAM solutions for $L-\lambda$. Moreover~: $\Nc^*_{[\lambda]}(H)=T_\lambda(\bigcup \tilde\Ic(u_-, u_+))$ where the union is taken on all the pairs of conjugate weak KAM solutions for $L-\lambda$.\\
 
 An immediate corollary of all these  results is the following~:  if $\pi^*(\Ac_{[\lambda]}^*(H))=M$, then there is a unique negative weak KAM solution $u$ and a unique positive weak KAM solution for $L-\lambda$, they are equal and $C^{1, 1}$ (i.e. $C^1$ with a Lipschitz derivative). In this case, we have~: $\Ac^*_{[\lambda ]}(H)= \Nc^*_{[\lambda]}(H)$ is the graph of $\lambda+du$.

\section{Proof of theorem \ref{T1}}
We assume that $H$ is a Tonelli Hamiltonian such that $\Nc^T_*(H)=T^*M$.\\
In order to prove theorem \ref{T1}, we begin by proving  that the periodic orbits are on some  invariant totally periodic Lagrangian graphs~:
\begin{prop}\label{P4}
For every closed $1$-form $\lambda$ of $M$, for every $(q_0,p_0)\in T^*M$ that is $T$-periodic for a certain $T>0$ and whose orbit  under the Hamiltonian flow is minimizing for $L-\lambda$, then $(q_0, p_0)$  belongs to a $C^1$ invariant Lagrangian graph $\Tc$ such that the orbit of every element of $\Tc$ is $T$-periodic, homotopic to the one of $(q_0, p_0)$  and has the same action for the Lagrangian $L-\lambda$ as the orbit of $(q,p)$. Moreover, $\Tc$ is the graph of  a closed 1-form that has the same cohomology class as $\lambda$.
\end{prop}
\demo
Let us consider $(q_0, p_0)$ as in the statement. Then, if we denote the cohomology class of $\lambda$ by $[\lambda$], we have~: $(q_0, p_0)\in \Nc^*_{[\lambda ]}(H)$, i.e. $(q_0, p_0)$ belongs to the Ma\~n\'e set associated to the cohomology class of $\lambda$. Let us use the notation~: $\gamma_0(t)=\pi\circ\varphi_t(q_0, p_0)$.

Because of Tonelli theorem, we know that for every $q\in M$, there exists a piece of orbit $(\varphi_t(q, p))_{t\in [0, T]}$ such that, if we denote the projection of this piece of orbit by $\gamma_q$ (i.e. $\gamma_q (t)=\pi\circ \varphi_t(q,p)$), then we have~:
\begin{enumerate}
\item[$\bullet$] $\gamma_q (T)=\gamma_q (0)=q$;
\item[$\bullet$] $\gamma_q$ is homotopic to $\gamma_0$;
\item[$\bullet$] for every absolutely continuous arc $\eta~: [0, T]\rightarrow M$ that is homotopic to $\gamma_0$ and such that~: $\eta (0)=\eta (T)=q$, we have~:
$\int_0^T(L(\gamma_q, \dot\gamma_q)-\lambda (\dot\gamma_q))\leq \int_0^T(L(\eta, \dot\eta)-\lambda (\dot\eta))$.
\end{enumerate}
  As every point on $T^*M$ is in some Ma\~n\'e set, then the orbit of every point has to be a graph by proposition \ref{Pr8}. We deduce that~: $\varphi_T(q, p)=(q,p)$, hence $(q, p)$ is a $T$-periodic point. It defines an invariant probability measure $\mu_q$, the one equidistributed along this orbit, defined by~:
  $$\forall f\in C^0(T^*M, \R), \int fd\mu=\frac{1}{T}\int_0^Tf\circ \varphi_t(q,p)dt.$$
  As the support of this measure is in some Ma\~n\'e set, this measure is minimizing for $L +\nu$ where $\nu$ is some closed 1-form. The rotation vector of this measure is   $\frac{1}{T}[\gamma]=\frac{1}{T}[\gamma_0]$ where we denote the homology class of $\gamma$ by $[\gamma]$; hence, having the same rotation vector, the supports of the measures $\mu_q$ and $\mu_{q_0}$ belong to the same Mather set and the support of $\mu_q$ is in $\Nc^*_{[\lambda ]}(H)$. We deduce that~:
  $$\forall q\in M, -T\alpha ([\lambda])=\int_0^T(L(\gamma, \dot\gamma)-\lambda( \dot\gamma))=\int_0^T(L(\gamma_0, \dot\gamma_0)-\lambda (\dot\gamma_0))$$
  because all these measures are minimizing for $L-\lambda$.
  
Finally,  for all  $q\in M$, we have found a point $(q,p)$ that is in the Mather set $\Mc^*_{[\lambda ]}(H)$.  As the Mather set is a Lipschitz graph, then the set of these points $(q,p)$ is a Lipschitz graph and coincides with the Mather set $\Mc^*_{[\lambda ]}(H)$. Moreover, we know that the Aubry set is a graph that contains the Mather set. Hence $\Ac^*_{[\lambda ]}(H)=\Mc^*_{[\lambda ]}(H)$. In this case, this set is the graph of a Lipschitz closed 1-form whose cohomology class is $[\lambda]$ (see subsection \ref{secweak}). As the dynamic restricted to this  $C^0$-Lagrangian graph is totally periodic, i.e. as $\varphi_{T|\Ac^*_{[\lambda ]}(H)}={\rm Id}_{\Ac^*_{[\lambda ]}(H)}$, we know   that this graph is in fact $C^1$ (this is proved in \cite{Arna3} by way of the so-called Green bundles).
 
\enddemo
We can apply this proposition to every periodic orbit. Indeed, such a periodic orbit is always contained in some Ma\~n\'e set $\Nc_c^*(H)$. We deduce from the previous proposition that $\Ac_c^*(H)$ is a $C^1$ Lagrangian graph, and that all the orbits contained in $\Ac_c^*(H)$ are periodic with the same period and are homotopic to each other. Moreover, we have seen in subsection \ref{secweak} that when the Aubry set is a graph above the whole zero section, then it coincides with the Ma\~n\'e set. Hence, we have proved that  $\Nc_c^*(H)$ is a $C^1$ Lagrangian graph, and that all the orbits contained in $\Nc_c^*(H)$ are periodic with the same period and are homotopic to each other.
\medskip

Let us know explain what happens to the other Ma\~n\'e sets, that correspond to the other cohomology classes.

\begin{prop} \label{P5} For every cohomology class $c\in H^1(M, \R)$, we have~: $\Ac_c^*(H)=\Nc_c^*(H)$ is the graph $\Gc_c$ of a Lipschitz closed 1-form.

\end{prop}
 
\demo
Let us assume that $(q, p)\in \Ac_c^*(H)$. Let $\lambda$ be a closed 1-form such that $[\lambda]=c$. Then there exists a sequence $(T_n)$ tending to $+\infty$  and a sequence $(\gamma_n)$ of absolutely continuous arcs $\gamma_n~: [0, T_n]\rightarrow M$ that are minimizing, such that $\gamma (0)=\gamma (T_n)=q$ and such that~:
$\displaystyle{\lim_{n\rightarrow \infty}\int_0^{t_n}(L(\gamma_n(t), \dot\gamma_n (t))-\lambda (\dot\gamma_n (t))+\alpha (c))dt=0}$ where $\alpha$ designates the $\alpha$ function of Mather. As every $\gamma_n$ is minimizing, it is the projection of a piece of orbit~: $\gamma_n(t)=\pi\circ\varphi_t(q, p_n)$. The corresponding orbit, being in a certain Ma\~n\'e set,  has to be a graph, hence it is periodic~: $\varphi_{t_n}(q, p_n)=(q, p_n)$. Moreover, we know (see proposition \ref{Pradial}) that in this case~: $\displaystyle{\lim_{n\rightarrow \infty} (q, p_n)=(q,p)}$.\\
We can use proposition \ref{P4}. Let $c_n\in H^1(M, \R)$ be the cohomology class such that $(q, p_n)\in \Nc_{c_n}^*(H)$. Then there exists a closed 1-form $\lambda_n$, whose cohomology class is $c_n$, so that $\Nc_{c_n}^*(H)$ is the graph of $\lambda_n$. We have in particular~: $p_n=\lambda_n(q)$ and $\displaystyle{p=\lim_{n\rightarrow \infty}\lambda_n(q)}$. Let us now prove that for every $Q\in M$, the sequence $(Q, \lambda_n(Q))$   converges to some point $(Q, P)$ that belongs to $\Ac_c^*(H)$. We will deduce that $\Ac_c^*(H)=\Nc_c^*(H)$ is the graph of a Lipschitz closed 1-form and then the proposition.

So let us consider $Q\in M$. For every $n\in\N$, we know by proposition \ref{P4} that $(Q, \lambda_n(Q))$ is $t_n$-periodic and that if we denote the projection of its orbit by $\Gamma_n(t)=\pi\circ \varphi_t(Q, \lambda_n(Q))$, then we have~:
\begin{enumerate}
\item[$\bullet$] $\Gamma_n$ is homotopic to $\gamma_n$;
\item[$\bullet$] $\int_0^{t_n}(L(\Gamma_n(t), \dot\Gamma_n(t))-\lambda_n(\dot\Gamma_n(t)))dt=
\int_0^{t_n}(L(\gamma_n(t), \dot\gamma_n(t))-\lambda_n(\dot\gamma_n(t)))dt$.
\end{enumerate}
We can then compute (the notation $[\lambda][\gamma]$ is just the usual product of a cohomology class with a homology class)~:\\
$\int_0^{t_n}(L(\Gamma_n(t), \dot\Gamma_n (t))-\lambda( \dot\Gamma_n (t))+\alpha (c))dt=$\\
$\int_0^{t_n}(L(\Gamma_n(t), \dot\Gamma_n (t))-\lambda_n ( \dot\Gamma_n (t)))dt -[\lambda-\lambda_n][\Gamma_n]+\alpha (c)t_n
=$\\
$\int_0^{t_n}(L(\gamma_n(t), \dot\gamma_n (t))-\lambda_n (\dot\gamma_n (t)))dt -[\lambda-\lambda_n][\gamma_n]+\alpha (c)t_n
=$\\
$\int_0^{t_n}(L(\gamma_n(t), \dot\gamma_n (t))-\lambda (\dot\gamma_n (t))+\alpha (c))dt$.\\
Then~:
$\displaystyle{\lim_{n\rightarrow \infty}\int_0^{t_n}(L(\Gamma_n(t), \dot\Gamma_n (t))-\lambda( \dot\Gamma_n (t))+\alpha (c))dt =}0$. By proposition \ref{Pradial}, this implies that $Q$ belongs to the projected Aubry set $\pi(\Ac_c^*(H))$ and that the sequence $(Q, \lambda_n(Q))$ converges to the unique point of $\Ac_c^*(H)$ that is above $Q$.

\enddemo
\begin{prop}
With the previous notations, the graphs $\Gc_c$ are disjoints~:
$$\forall c, d\in H^1(M, \R), c\not=d\Rightarrow \Gc_c\cap \Gc_d=\emptyset.$$
\end{prop}
\demo
We borrow the main elements of the proof to \cite{massart1}. Let us assume that there exists $c, d\in H^1(M, \R)$ such that $\Gc_c\cap \Gc_d\not= \emptyset$. Then $\Gc_c\cap\Gc_d$ is a compact invariant subset  and there exists an invariant Borel probability measure $\mu^*$ (dual of $\mu$) whose support is contained in $\Gc_c\cap \Gc_d$.  Hence $\mu$ is minimizing for $L-\lambda $ and $L-\eta$ if $[\lambda]=c$ and $[\eta]=d$~:
$$\int(L-\lambda+\alpha(c))d\mu=0\quad{\rm and}\quad \int(L-\eta+\alpha (d))d\mu=0.$$
We deduce that  for every $t\in[0, 1]$, we have~:
$$\int(L-(t\lambda +(1-t)\eta)+t\alpha (c)+(1-t)\alpha (d))d\mu=0$$
and then~: $\alpha (tc+(1-t)d)\geq -\int(L-(t\lambda+(1-t)\eta))d\mu =t\alpha (c)+(1-t)\alpha (d)$. As the function $\alpha$ is convex, this implies~: $\alpha (tc+(1-t)d)=t\alpha (c)+(1-t)\alpha (d)$. Hence $\mu$ is minimizing for $L-(t\lambda +(1-t)\eta)$. This implies that the support of $\mu^*$ is contained in $\Mc_{tc+(1-d)}^*(H)\subset\Ac_{tc+(1-d)}^*(H)=\Nc_{tc+(1-d)}^*(H)=\Gc_{tc+(1-d)}$. Let us now consider $(q, p)\in \Gc_{\frac{1}{2}(c+d)}$. As $(q, p)$ belongs to $ \Ac_{\frac{1}{2}(c+d)}^*(H)$, there exists a sequence $(T_n)$ tending to $+\infty$ and a sequence of $C^1$ arcs $\gamma_n~: [0, T_n]\rightarrow M$ such that $\gamma_n(0)=\gamma_n(T_n)=q$ and~:
$$\lim_{n\rightarrow \infty}(\int_0^{T_n}(L(\gamma_n(t), \dot\gamma_n(t))-\frac{1}{2}(\lambda +\eta)(\dot\gamma_n(t))+\alpha (\frac{1}{2}(c+d))dt)=0.$$
The left term of the previous equality is the limit of the sum of two terms~: \\
$\frac{1}{2}\int_0^{T_n}(L(\gamma_n(t),\dot\gamma_n(t))-\lambda (\dot\gamma_n(t))+\alpha (c))dt$ and $\frac{1}{2}\int_0^{T_n}(L(\gamma_n(t),\dot\gamma_n(t))-\eta (\dot\gamma_n(t))+\alpha (d))dt$, each of these terms being non negative. We deduce that~:\begin{enumerate}
 \item[$\bullet$] $\displaystyle{\lim_{n\rightarrow \infty}\int_0^{T_n}(L(\gamma_n(t),\dot\gamma_n(t))-\lambda (\dot\gamma_n(t))+\alpha (c))dt=0}$; 
 \item[$\bullet$] $\displaystyle{\lim_{n\rightarrow \infty}\int_0^{T_n}(L(\gamma_n(t),\dot\gamma_n(t))-\eta (\dot\gamma_n(t))+\alpha (d))dt=0}$;
\end{enumerate}
and by proposition \ref{Pradial}~: $$\lim_{n\rightarrow \infty}(\gamma_n(0), \dot\gamma_n(0))\in \Ac_c(H)\cap \Ac_d(H).$$
We have finally proved that $ \Gc_{\frac{1}{2}(c+d)}=\Ac_{\frac{1}{2}(c+d)}^*(H)\subset \Ac_c^*(H)\cap \Ac^*_d(H)=\Gc_c\cap\Gc_d$, hence the two graphs $\Gc_c$ and $\Gc_d$ are equal, and their cohomology classes are also equal~: $c=d$. 
\enddemo
Let us now finish the proof of theorem \ref{T1}. We have found a partition  of $T^*M$ into Lipschitz Lagrangian graphs $(\Gc_c)_{c\in H^1(M, \R))}$, where $\Gc_c$ is the graph of a Lipschitz 1-form whose cohomology class is $c$ and is  equal to $\Ac_c^*(H)=\Nc_c^*(H)$. Each Ma\~n\'e set being chain recurrent, we deduce that the dynamic restricted to each $\Gc_c$ is chain recurrent. \\
We are then exactly in the case of a $C^0$-integrable Hamiltonian that we described in \cite{Arna3}. We can   apply the results of \cite{Arna3} and deduce that there exists a dense $G_\delta$-subset of $T^*M$ filled by invariant $C^1$ Lagrangian graphs. Finally, let us notice that it is proved in \cite{masor} that the $\beta$ function of every $C^0$ integrable Tonelli Hamiltonian is differentiable everywhere.This ends the proof of the implication~: $2\Rightarrow 1$.

Let us now prove that $1\Rightarrow 2$. We assume that there is a partition of $T^*M$ into invariant Lagrangian Lipschitz graph. Then to each of these Lipschitz graphs corresponds a $C^{1, 1}$ weak KAM solution and then the orbit of every point of this graph is in some Ma\~n\'e set. This implies~: $T^*M=\Nc^T_*(M)$.

 \section{Proof of the corollaries}
\subsection{Proof of corollary \ref{C2}} We only have to prove that $2\Rightarrow 1$. We assume that $T^*M$ is covered by the union of the invariant Lipschitz Lagrangian graphs (resp. smooth Lagrangian graphs). Then to each of these Lipschitz graphs correspond a $C^{1, 1}$ weak KAM solution and then the orbit of every point of this graph is in some Ma\~n\'e set. This implies~: $T^*M=\Nc^T_*(H)$. We can apply theorem \ref{T1} and proposition \ref{P5}. Then there exists a partition of $T^*M$ into Lipschitz Lagrangian graphs $(\Gc_c)_{c\in H^1(M, \R))}$, where $\Gc_c$ is the graph of a Lipschitz 1-form whose cohomology class is $c$ and is  equal to $\Ac_c^*(H)=\Nc_c^*(H)$. Let us look at what happens in the smooth case~: if $N$ is one of the smooth invariant Lagrangian graphs, then it is contained in some Ma\~n\'e set and then is equal to some $\Gc_c$. We obtain then that there is a partition of $T^*M$ into some smooth  $\Gc_c$. As $(\Gc_c)_{c\in H^1(M, \R)}$ is a partition of $T^*M$, we deduce that all the $\Gc_c$ are smooth.

\subsection{Proof of corollary \ref{C3}} We just have to prove that $2\Rightarrow 1$. We assume $T^*M$ is covered by the union of its Lagrangian invariant smooth submanifolds that are Hamiltonianly isotopic to some smooth Lagrangian graph. We have proved in \cite{Arna2} a multidimensional Birkhoff theorem~: every Lagrangian invariant smooth submanifold that is  Hamiltonianly isotopic to some smooth Lagrangian graph is a smooth graph. Then corollary \ref{C3} becomes a corollary of corollary \ref{C2}.

\end{document}